\documentclass[12pt]{amsart}

\usepackage[utf8]{inputenc}
\usepackage[a4paper,nomarginpar]{geometry}
\usepackage[english]{babel}
\usepackage{enumitem}
\usepackage{amsmath,amsfonts,amssymb,amsthm,amscd}
\usepackage{tikz}
\usepackage{xcolor}

\allowdisplaybreaks

\theoremstyle{theorem}
\newtheorem{theorem}{Theorem}%[section]

\theoremstyle{definition}
\newtheorem{definition}[theorem]{Definition}

\newtheorem{remark}[theorem]{Remark}

\theoremstyle{remark} \theoremstyle{question} \theoremstyle{example}

\newcommand{\N}{\mathbb{N}}
\newcommand{\Z}{\mathbb{Z}}

\newcommand{\R}{\mathbb{R}}
\newcommand{\C}{\mathbb{C}}
\newcommand{\D}{\mathbb{D}}

\newcommand{\K}{\mathbb{K}}

\newcommand{\cF}{\mathcal{F}}
\newcommand{\cG}{\mathcal{G}}

\newcommand{\cL}{\mathcal{L}}

\newcommand{\cN}{\mathcal{N}}

\newcommand{\vertiii}[1]{{\left\vert\kern-0.25ex\left\vert\kern-0.25ex\left\vert #1 
    \right\vert\kern-0.25ex\right\vert\kern-0.25ex\right\vert}}

\newcommand{\Lip}{\operatorname{Lip}}

\newcommand{\eps}{\varepsilon}
\newcommand{\ov} {\overline}
\newcommand{\wt}{\widetilde}

\begin{document}


\title[ON THE STABILITY OF GENERALIZED HYPERBOLIC OPERATORS]{ON THE STABILITY OF\\ GENERALIZED HYPERBOLIC OPERATORS}

\subjclass[2020]{Primary 47A99, 37C20, 37B25; Secondary 47A16. }
\keywords{Hyperbolicity; Generalized hyperbolicity; Structural stability; Time-dependent stability; Linear operators; Banach spaces.}
\date{}
\dedicatory{}
\maketitle

\begin{center}
{\sc Nilson C. Bernardes Jr.}

\medskip
Institut Universitari de Matem\`atica Pura i Aplicada\\
Universitat Polit\`ecnica de Val\`encia\\
Cam\'i de Vera S/N, Edifici 8E, Acces F, 4a Planta, Val\`encia, 46022, Spain\\
and\\
Departamento de Matem\'atica Aplicada, Instituto de Matem\'atica\\
Universidade Federal do Rio de Janeiro, Caixa Postal 68530\\
Rio de Janeiro, RJ, 21941-909, Brazil

\smallskip
{\it e-mail}: ncbernardesjr@gmail.com

\end{center}

\medskip

\begin{abstract}
Our main goal is to prove that every invertible generalized hyperbolic operator on a Banach space has a stability property,
known as time-dependent stability, which was introduced by J. M. Franks (Invent.\ Math.\ {\bf 24} (1974), 163--172) and is
stronger than structural stability.
\end{abstract}

%%%%%%%%%%%%%%%%%%%%%%%%%%%%%%%%%%%%%%%%%%%%%%%%%%%%%%%%%%%%%%%
%%%%%%%%%%%%%%%%%%%%%%%%%%%%%%%%%%%%%%%%%%%%%%%%%%%%%%%%%%%%%%%

\section{Introduction}

Structural stability is one of the most important concepts in the modern theory of dynamical systems and differential equations.
It was originally introduced by A. Andronov and L. Pontrjagin \cite{AAndLPon37} for a certain class of differentiable flows on the plane.
Later, variations of the original concept were introduced and investigated in different contexts.
In the area of differentiable dynamics, a diffeomorphism $f : M \to M$ on a compact manifold is said to be {\em structurally stable} 
if every diffeomorphism $g : M \to M$ that is sufficiently close to $f$ in the $C^1$-topology is {\em topologically conjugate} to $f$, that is,
there exists a homeomorphism $h : M \to M$ such that $g = h \circ f \circ h^{-1}$.
The importance of this concept lies in the fact that models of physical systems are only approximations of the true systems
and it is important to know whether the qualitative behaviors of the approximation and the true dynamical system are the same.
We refer the reader to \cite{AKatBHas95,JRob72,MShu87} for nice expositions about structural stability.

J. Franks \cite{JFra74} observed that, in most physical situations, it seems likely that the true dynamical system is not really autonomous
but is instead to a certain extent time dependent.
So, he proposed the following concept of stability:
A diffeomorphism $f : M \to M$ on a compact manifold is said to be {\em time-dependent stable} if there is a neighborhood $\cN$ of $f$ 
in the $C^1$-topology such that $g_1 \circ \cdots \circ g_p$ is topologically conjugate to $f^p$ whenever $g_1,\ldots,g_p \in \cN$ 
(an important point here is that $\cN$ is independent of $p$).
It is clear that this concept is stronger than structural stability.

In the present note we are interested in the dynamics of invertible operators on Banach spaces.
Our main motivation comes from the following classical result.

\medskip\noindent
{\bf Hartman's theorem.} {\it Every invertible hyperbolic operator on a Banach space is structurally stable.}

\medskip
The name ``Hartman's theorem'' for the above result was coined by C. Pugh \cite{CPug69}.
Actually, P. Hartman \cite{PHar60} established it for operators on finite-dimensional euclidean spaces.
The extension to arbitrary Banach spaces was independently obtained by J. Palis \cite{JPal68} and C. Pugh \cite{CPug69}
(both were motivated by an argument due to J. Moser \cite{JMos69}).
This theorem is the main tool for proving the celebrated {\em Grobman-Hartman theorem} 
on the local linearization of diffeomorphisms around hyperbolic fixed points in Banach spaces \cite{JPal68,CPug69}.
We refer the reader to \cite{JGucPHol83,AKatBHas95,WMelJPal82,MShu87} for more details on this important
{\em linearization theorem} and its applications.

A basic question is whether or not the converse of {\em Hartman's theorem} is always true.
It was soon realized that the answer is positive in the finite-dimensional setting (see J. Robbin \cite{JRob72}),
but the full question remained open for more than 50 years.
This question was finally answered in the negative in the recent paper \cite{NBerAMes21} of the author with A. Messaoudi,
where the first examples of structurally stable operators that are not hyperbolic were obtained.
A little later, the same authors proved that every invertible generalized hyperbolic operator on a Banach space is structurally stable
\cite[Theorem~1]{NBerAMes20}, which implied a {\em generalized Grobman-Hartman theorem} \cite[Theorem~3]{NBerAMes20}.
Moreover, in the recent papers \cite{BerCarDarFavPer25,KLeeCMor24}, it was independently shown that 
every invertible generalized hyperbolic operator on a Banach space is also {\em topologically stable},
a well-known concept of stability introduced by P. Walters~\cite{PWal70}.

In short, the significant developments in the field of linear dynamics in recent decades have motivated some authors 
to revisit stability theory in this context \cite{FBay21,BerCarDarFavPer25,NBerAMes20,NBerAMes21,KLeeCMor24,FLenAMesTA}. 
Several relevant results have already been obtained, and we believe that this branch of linear dynamics will continue 
to undergo significant development for many years to come. 
The present note contains another contribution to this line of investigation.
We address the following basic question related to {\em Hartman's theorem}: 
\begin{quote}
{\it Is every invertible hyperbolic operator on a Banach space time-dependent stable?}
\end{quote}

Our main goal is to give a positive answer to this question. Actually, we will establish the following more general result.

\begin{theorem}\label{Main}
Every invertible generalized hyperbolic operator on a Banach space is time-dependent stable.
\end{theorem}

In fact, this theorem will be derived from an even more general result (Theorem~\ref{Main2}).

The note is organized as follows.
In Section~\ref{Section2} we present all the necessary preliminaries for the present work.
In Section~\ref{Section3} we prove a refinement of Theorem~\ref{Main}, namely Theorem~\ref{Main2}, 
in which we also analyze the problem of uniqueness of the homeomorphism that establishes the conjugation
and the possibility of choosing this homeomorphism in such a way that we have continuous dependence with respect to the perturbation.
Some additional remarks are also presented at the end of the section.

%%%%%%%%%%%%%%%%%%%%%%%%%%%%%%%%%%%%%%%%%%%%%%%%%%%%%%%%%%%%%%%
%%%%%%%%%%%%%%%%%%%%%%%%%%%%%%%%%%%%%%%%%%%%%%%%%%%%%%%%%%%%%%%

\section{Preliminaries}\label{Section2}

Throughout $\K$ denotes either the field $\R$ of real numbers or the field $\C$ of complex numbers,
$\N$ denotes the set of all positive integers and $\N_0\!:= \N \cup \{0\}$.

All Banach spaces are assumed to be over $\K$, unless otherwise specified.

Given a continuous linear map $T$ from a Banach space $X$ into a Banach space $Y$, 
recall that the {\em norm} of $T$ is the non-negative real number given by
\[
\|T\| = \sup_{\|x\| \leq 1} \|T(x)\|.
\]
This notion will play an important role in the next section.
If $F : X \to Y$ is any map, we define
\[
\|F\|_\infty = \sup_{x \in X} \|F(x)\| \ \ \text{ and } \ \ \Lip(F) = \sup_{x \neq x'} \frac{\|F(x) - F(x')\|}{\|x - x'\|}\,\cdot
\]
We say that $F$ is a {\em bounded map} if $\|F\|_\infty < \infty$ and it is a {\em Lipschitz map} if $\Lip(F) < \infty$.

A continuous linear map from a Banach space $X$ into itself will be called simply an {\em operator} on $X$.
We denote by $I_X$ the {\em identity operator} on $X$.
In this note we are interested in the dynamics of {\em invertible operators}.

Recall that an invertible operator $T$ on a Banach space $X$ is said to be {\em hyperbolic} if its spectrum 
$\sigma(T)$ does not intersect the unit circle in the complex plane. 
It is well known that $T$ is hyperbolic if and only if there are an equivalent norm $\vertiii{\cdot}$ on $X$ and a direct sum decomposition
\[
X = X_s \oplus X_u, \ T = T_s \oplus T_u, 
\]
called the {\em hyperbolic splitting} of $T$, such that $X_s$ and $X_u$ are closed $T$-invariant subspaces of $X$ 
(the {\em stable} and the {\em unstable subspaces} for $T$), 
$T_s = T|_{X_s}$ is a proper contraction (i.e., $\|T_s\| < 1$) and $T_u = T|_{X_u}$ is a proper dilation (i.e., $\|T_u^{-1}\| < 1$),
where the norms $\|T_s\|$ and $\|T_u^{-1}\|$ are calculated with respect to the new norm $\vertiii{\cdot}$ on $X$.

Recall that an invertible operator $T$ on a Banach space $X$ is said to be {\em generalized hyperbolic} 
if there is a direct sum decomposition
\[
X = M \oplus N,
\]
where $M$ and $N$ are closed subspaces of $X$ with the following properties:
\begin{itemize}
\item [(GH1)] $T(M) \subset M$ and $\sigma(T|_M) \subset \D$,
\item [(GH2)] $T^{-1}(N) \subset N$ and $\sigma(T^{-1}|_N) \subset \D$,
\end{itemize}
where $\D$ denotes the open unit disc in the complex plane.
This class of operators appeared in the statement of \cite[Theorem~A]{BerCirDarMesPuj18}, which was the major tool to prove 
the existence of operators that have the shadowing property but are not hyperbolic \cite[Theorem~B]{BerCirDarMesPuj18}.
However, the terminology ``generalized hyperbolic'' was introduced in \cite{CirGolPuj21},
where additional dynamical properties of these operators were investigated 
(see also \cite{BerCarDarFavPer25,NBerAMes20,NBerAPer24}).

In the definition of the concept of structural stability, it is usual to consider perturbations that are small in the sense of the $C^1$-topology. 
However, in certain situations, some authors allow more general perturbations, namely: 
Lipschitz perturbations with small sup norm and small Lipschitz constant (see \cite{CPug69}, for instance).
As in \cite{NBerAMes20,NBerAMes21}, we adopt the following definition of structural stability in the setting of linear dynamics.

\begin{definition}
An invertible operator $T$ on a Banach space $X$ is said to be {\em structurally stable} if there exists $\eps > 0$ such that 
\[
T + L \text{ is topologically conjugate to } T
\]
whenever $L : X \to X$ is a Lipschitz map with $\|L\|_\infty < \eps$ and $\Lip(L) < \eps$.
\end{definition}

Therefore, it is natural to adopt the following definition of time-dependent stability in the present context.

\begin{definition}
An invertible operator $T$ on a Banach space $X$ is said to be {\em time-dependent stable} if there exists $\eps > 0$ such that 
\[
(T+L_1) \circ \cdots \circ (T+L_p) \text{ is topologically conjugate to } T^p
\]
whenever $p \in \N$ and $L_1,\ldots,L_p$ are Lipschitz maps from $X$ into itself with 
$\|L_j\|_\infty < \eps$ and $\Lip(L_j) < \eps$ for all $j \in \{1,\ldots,p\}$.
\end{definition}

%%%%%%%%%%%%%%%%%%%%%%%%%%%%%%%%%%%%%%%%%%%%%%%%%%%%%%%%%%%%%%%
%%%%%%%%%%%%%%%%%%%%%%%%%%%%%%%%%%%%%%%%%%%%%%%%%%%%%%%%%%%%%%%

\section{Main result}\label{Section3}

Let us now state and prove our main result, which clearly implies Theorem~\ref{Main}.

\begin{theorem}\label{Main2}
For any invertible generalized hyperbolic operator $T$ on any Banach space $X$, the following properties hold:
\begin{itemize}
\item [(A)] There exists $\eps > 0$ such that for any $p \in \N$ and any bounded Lipschitz maps $L_j : X \to X$, $j = 1,\ldots,p$, with
\[
\Lip(L_j) < \eps \ \text{ for all } j \in \{1,\ldots,p\},
\]
there exists a homeomorphism $h : X \to X$ satisfying
\[
(T+L_1) \circ \cdots \circ (T+L_p) = h \circ T^p \circ h^{-1} \ \text{ and } \ \ \|h - I_X\|_\infty < \infty.
\]
\item [(B)] Given $\delta > 0$, there exists $\eps > 0$ such that for any $p \in \N$ and any Lipschitz maps $L_j : X \to X$, 
$j = 1,\ldots,p$, with
\[
\|L_j\|_\infty < \eps \ \text{ and } \ \Lip(L_j) < \eps \ \text{ for all } j \in \{1,\ldots,p\},
\]
there exists a homeomorphism $h : X \to X$ satisfying
\[
(T+L_1) \circ \cdots \circ (T+L_p) = h \circ T^p \circ h^{-1} \ \text{ and } \ \ \|h - I_X\|_\infty < \delta.
\]
\end{itemize}
Moreover, the homeomorphism $h$ given in (A) (or (B)) is unique if and only if the operator $T$ is hyperbolic.
Nevertheless, it is always possible to choose such a homeomorphism $h_\cL$ for each $p$-tuple $\cL = (L_1,\ldots,L_p)$ of maps
as in (A) (resp.\ (B)) in such a way that the correspondence $\cL \to h_\cL$ is uniformly continuous in the following sense:
given $\gamma > 0$, there exists $\eta > 0$ such that 
\[
\|h_\cL - h_{\cL'}\|_\infty < \gamma
\]
whenever $\cL = (L_1,\ldots,L_p)$ and $\cL' = (L'_1,\ldots,L'_p)$ are $p$-tuples of maps as in (A) (resp.\ (B)) with
\[
\|L_j - L'_j\|_\infty < \eta \ \text{ for all } j \in \{1,\ldots,p\}.
\]
\end{theorem}

Note that the assumption that the $L_j$'s are small in the sup norm is not necessary to guarantee that
$(T+L_1) \circ \cdots \circ (T+L_p)$ is topologically conjugate to $T^p$,
the important assumption is that the $L_j$'s are bounded maps with small Lipschitz constants.
However, by adding the hypothesis that the $L_j$'s are small in the sup norm,
we can ensure that the homeomorphism that establishes the conjugation can be chosen to be uniformly close to the identity map.

In the proof below we combine ideas from \cite{NBerAMes20} and \cite{JFra74} with some new ideas.

\begin{proof}
Fix an invertible generalized hyperbolic operator $T$ on a Banach space $X$. Let 
\begin{equation}\label{eq1}
X = M \oplus N
\end{equation}
be a direct sum decomposition given by the definition of generalized hyperbolicity and consider the closed subspace of $X$ given by
\[
Y = M + T^{-1}(N).
\]
By the spectral radius formula, there exist $a \geq 1$ and $t \in (0,1)$ such that
\begin{equation}\label{eq2}
\|T^n(y)\| \leq a\,t^n\|y\| \ \text{ and } \ \|T^{-n}(z)\| \leq a\,t^n\|z\|
\ \text{ whenever } n \in \N_0, y \in M, z \in N.
\end{equation}
Let $P_M : X \to M$ and $P_N : X \to N$ denote the canonical projections associated to the direct sum decomposition (\ref{eq1}) and define 
\begin{equation}\label{eq3}
b = \max\{\|P_M\|,\|P_N\|\}.
\end{equation}
Choose $\delta \in (0,1)$. We shall prove that the number
\[
\eps = \min\Big\{\frac{1-t}{a\,b\,(1+t)},\frac{1}{\|T^{-1}\|}\Big\}\, \delta
\]
satisfies properties (A) and (B). Indeed, fix $p \in \N$ and a $p$-tuple $\cL = (L_0,\ldots,L_{p-1})$ of bounded Lipschitz maps 
$L_j : X \to X$ such that
\[
\Lip(L_j) < \eps \ \text{ for all } j \in \{0,\ldots,p-1\}.
\]
Let $S_j = T + L_j$ for each $j$. Since $\Lip(L_j) < \|T^{-1}\|^{-1}$, we have that $S_j$ is a homeomorphism.
We shall prove that $S_{p-1} \circ \cdots \circ S_0$ is topologically conjugate to $T^p$.
For this purpose, let $\Z_p = \{0,\ldots,p-1\}$ endowed with the operation $+_p$ of addition mod $p$ and with the trivial metric $d_0$.
Consider the product $\wt{X} = X \times \Z_p$ endowed with the metric
\[
d((x,j),(x',j')) = \max\{\|x - x'\|,d_0(j,j')\}.
\]
For $(x,j),(x',j) \in \wt{X}$ and $\lambda \in \K$, we define
\[
(x,j) \pm (x',j) = (x \pm x',j), \ \ \ \lambda (x,j) = (\lambda x,j) \ \ \text{ and } \ \ |(x,j)| = \|x\|.
\]
Note that the above operations of ``addition'' and ``subtraction'' on $\wt{X}$ are defined only for pairs of points with the same second coordinate.
Moreover, we define
\[
\sum_{k = 0}^\infty (x_k,j) = \Big(\sum_{k = 0}^\infty x_k\,,j\Big),
\]
provided $\sum_{k = 0}^\infty x_k$ is a convergent series in $X$. We will need the function spaces defined below:
\begin{itemize}
\item $\cF$ is the set of all continuous bounded maps $F : \wt{X} \to \wt{X}$ such that
  \[
  F(X \times \{j\}) \subset Y \times \{j\} \ \text{ for all } j \in \Z_p.
  \]
\item $\cG$ is the set of all continuous bounded maps $F : \wt{X} \to \wt{X}$ such that
  \[
  F(X \times \{j\}) \subset X \times \{j +_p 1\} \ \text{ for all } j \in \Z_p.
  \]
\end{itemize}
We recall that $F : \wt{X} \to \wt{X}$ is a {\em bounded map} if its range $F(\wt{X})$ is a bounded set in $\wt{X}$,
i.e., a set with finite diameter.
Although $\wt{X}$ is not a vector space and $|\cdot|$ is not a norm, both $\cF$ and $\cG$ are Banach spaces
when endowed with the algebraic operations defined pointwise
\[
(F_1 + F_2)(x,j) = F_1(x,j) + F_2(x,j) \ \ \text{ and } \ \ (\lambda F)(x,j) = \lambda F(x,j),
\]
and the supremum norm 
\[
|F|_\infty = \sup_{(x,j) \in \wt{X}} |F(x,j)|.
\]
Note that the zero vector of $\cF$ (resp.\ $\cG$) is the map given by
\[
O_\cF(x,j) = (0,j) \ \ \ \ (\text{resp.} \ O_\cG(x,j) = (0,j+_p 1)).
\]
Let $Q_M : \wt{X} \to M \times \Z_p$ and $Q_N : \wt{X} \to N \times \Z_p$ be given by
\[
Q_M(x,j) = (P_M(x),j) \ \ \text{ and } \ \ Q_N(x,j) = (P_N(x),j).
\]
Consider the homeomorphisms $\wt{T} : \wt{X} \to \wt{X}$ and $S : \wt{X} \to \wt{X}$ defined by
\[
\wt{T}(x,j) = (T(x),j +_p 1) \ \ \text{ and } \ \ S(x,j) = (S_j(x),j +_p 1).
\]
Given a homeomorphism $R : \wt{X} \to \wt{X}$ with $R(X \times \{j\}) = X \times \{j +_p 1\}$ for all $j$, 
we consider the map
\begin{equation}\label{eq4}
\Psi : F \in \cF \mapsto F \circ R - \wt{T} \circ F \in \cG.
\end{equation}
It is easy to check that $\Psi$ is linear and continuous. 

We claim that $\Psi$ is injective. Indeed, let $F \in \cF$ with $\Psi(F) = O_\cG$. We have to show that $F = O_\cF$.
The equality $\Psi(F) = O_\cG$ gives
\begin{equation}\label{E}
F \circ R^n = \wt{T}^n \circ F \ \ \text{ for all } n \in \N.
\end{equation}
Hence, for every $(x,j) \in \wt{X}$ and $n \in \N$,
\begin{equation}\label{A1}
F(x,j) = \wt{T}^n(F(R^{-n}(x,j))) = \alpha_n(x,j) + \beta_n(x,j),
\end{equation}
where
\[
\alpha_n(x,j) = \wt{T}^n(Q_M(F(R^{-n}(x,j)))) \ \text{ and } \ \beta_n(x,j) = \wt{T}^n(Q_N(F(R^{-n}(x,j)))).
\]
It is clear that, for each $n \in \N$,
\begin{equation}\label{A2}
\alpha_n(x,j) \in M \times \Z_p \ \ \text{ for all } (x,j) \in \wt{X}.
\end{equation}
Let us show that, for each $n \in \N$,
\begin{equation}\label{A3}
\beta_n(x,j) \in T^{-1}(N) \times \Z_p \ \ \text{ for all } (x,j) \in \wt{X}.
\end{equation}
Since $F(R^{-1}(x,j)) \in Y \times \Z_p$, we must have $Q_N(F(R^{-1}(x,j))) \in T^{-1}(N) \times \Z_p$, and so
\[
\beta_1(x,j) = \wt{T}(Q_N(F(R^{-1}(x,j)))) \in N \times \Z_p.
\]
Hence, since $\alpha_1(x,j) \in M \times \Z_p$ and $\alpha_1(x,j) + \beta_1(x,j) = F(x,j) \in Y \times \Z_p$, 
we see that (\ref{A3}) holds for $n = 1$. Assume that (\ref{A3}) holds for a certain $n \in \N$. Then,
\[
\beta_{n+1}(x,j) = \wt{T}(\beta_n(R^{-1}(x,j))) \in N \times \Z_p.
\]
As before, we conclude that (\ref{A3}) holds with $n+1$ in the place of $n$. By induction, (\ref{A3}) holds for all $n \in \N$. 
By (\ref{A1}), (\ref{A2}) and (\ref{A3}),
\[
Q_M(F(x,j)) = \alpha_n(x,j) = \wt{T}^n(Q_M(F(R^{-n}(x,j)))) \ \ \text{ for all } n \in \N.
\]
Thus, by (\ref{eq2}) and (\ref{eq3}),
\[
|Q_M(F(x,j))| \leq a\,b\,t^n |F|_\infty \ \ \text{ for all } n \in \N,
\]
implying that
\begin{equation}\label{A4}
Q_M(F(x,j)) = (0,j) \ \ \text{ for all } (x,j) \in \wt{X}.
\end{equation}
Now, by (\ref{E}), for every $(x,j) \in \wt{X}$ and $n \in \N$,
\begin{equation}\label{B1}
F(x,j) = \wt{T}^{-n}(F(R^n(x,j))) = \varphi_n(x,j) + \psi_n(x,j),
\end{equation}
where
\[
\varphi_n(x,j) = \wt{T}^{-n}(Q_M(F(R^n(x,j)))) \ \text{ and } \ \psi_n(x,j) = \wt{T}^{-n}(Q_N(F(R^n(x,j)))).
\]
It is clear that, for each $n \in \N$,
\begin{equation}\label{B2}
\psi_n(x,j) \in T^{-1}(N) \times \Z_p \ \ \text{ for all } (x,j) \in \wt{X}.
\end{equation}
Let us show that, for each $n \in \N$,
\begin{equation}\label{B3}
\varphi_n(x,j) \in M \times \Z_p \ \ \text{ for all } (x,j) \in \wt{X}.
\end{equation}
Since
\begin{align*}
\wt{T}(Q_N(\varphi_1(x,j))) 
&= \wt{T}(\varphi_1(x,j)) - \wt{T}(Q_M(\varphi_1(x,j)))\\
&= Q_M(F(R(x,j))) - \wt{T}(Q_M(\varphi_1(x,j))) \in M \times \Z_p,
\end{align*}
we have that
\begin{equation}\label{B4}
Q_N(\varphi_1(x,j)) \in T^{-1}(M) \times \Z_p.
\end{equation}
By (\ref{B1}) with $n = 1$, we can write
\[
F(x,j) = Q_M(\varphi_1(x,j)) + (Q_N(\varphi_1(x,j)) + \psi_1(x,j)).
\]
Since $Q_M(\varphi_1(x,j)) \in M \times \Z_p$ and $Q_N(\varphi_1(x,j)) + \psi_1(x,j) \in N \times \Z_p$, 
the fact that $F(x,j) \in Y \times \Z_p$ implies that $Q_N(\varphi_1(x,j)) + \psi_1(x,j) \in T^{-1}(N) \times \Z_p$.
Hence, by (\ref{B2}),
\begin{equation}\label{B5}
Q_N(\varphi_1(x,j)) \in T^{-1}(N) \times \Z_p.
\end{equation}
By (\ref{B4}) and (\ref{B5}), $Q_N(\varphi_1(x,j)) \in \{0\} \times \Z_p$, which gives (\ref{B3}) for $n = 1$.
Assume that (\ref{B3}) holds for a certain $n \in \N$. Then,
\begin{align*}
\wt{T}(Q_N(\varphi_{n+1}(x,j))) 
&= \wt{T}(\varphi_{n+1}(x,j)) - \wt{T}(Q_M(\varphi_{n+1}(x,j)))\\
&= \varphi_n(R(x,j)) - \wt{T}(Q_M(\varphi_{n+1}(x,j))) \in M \times \Z_p,
\end{align*}
and so
\[
Q_N(\varphi_{n+1}(x,j)) \in T^{-1}(M) \times \Z_p.
\]
By arguing as above, we obtain $Q_N(\varphi_{n+1}(x,j)) \in \{0\} \times \Z_p$, and so (\ref{B3}) holds with $n+1$ in the place of $n$.
By induction, (\ref{B3}) holds for all $n \in \N$. 
By (\ref{B1}), (\ref{B2}) and (\ref{B3}),
\[
Q_N(F(x,j)) = \psi_n(x,j) = \wt{T}^{-n}(Q_N(F(R^n(x,j)))) \ \ \text{ for all } n \in \N.
\]
Therefore, by (\ref{eq2}) and (\ref{eq3}),
\[
|Q_N(F(x,j))| \leq a\,b\,t^n |F|_\infty \ \ \text{ for all } n \in \N,
\]
which implies that
\begin{equation}\label{B6}
Q_N(F(x,j)) = (0,j) \ \ \text{ for all } (x,j) \in \wt{X}.
\end{equation}
By (\ref{A4}) and (\ref{B6}), we conclude that $F = O_\cF$, as desired. 

Let us now prove that $\Psi$ is surjective (hence bijective) and its inverse is given by
\begin{equation}\label{eq5}
\Psi^{-1}(G)(x,j) = \sum_{k=0}^\infty \wt{T}^k(Q_M(G(R^{-k-1}(x,j))))
                     - \sum_{k=1}^\infty \wt{T}^{-k}(Q_N(G(R^{k-1}(x,j)))).
\end{equation}
For this purpose, take $G \in \cG$. Note that the partial sums of both series in (\ref{eq5}) are well defined. 
Moreover, by (\ref{eq2}) and~(\ref{eq3}),
\begin{equation}\label{eq6}
\sum_{k=0}^\infty |\wt{T}^k(Q_M(G(R^{-k-1}(x,j))))| 
  \leq \sum_{k=0}^\infty a\, b\, t^k |G(R^{-k-1}(x,j))|
  \leq \frac{a\, b}{1-t}\, |G|_\infty
\end{equation}
and
\begin{equation}\label{eq7}
\sum_{k=1}^\infty |\wt{T}^{-k}(Q_N(G(R^{k-1}(x,j))))| 
  \leq \sum_{k=1}^\infty a\, b\, t^k |G(R^{k-1}(x,j))|
  \leq \frac{a\, b\, t}{1-t}\, |G|_\infty.
\end{equation}
This implies that the first (resp.\ second) series in (\ref{eq5}) is well defined and determines an element of $M \times \{j\}$ 
(resp.\ $T^{-1}(N) \times \{j\}$). Hence, the difference in (\ref{eq5}) is an element of $Y \times \{j\}$.
Since the convergence is uniform, we conclude that the right hand side of (\ref{eq5}) defines a map $F \in \cF$.
Now,
\begin{align*}
\Psi(F)(x,j) &= F(R(x,j)) - \wt{T}(F(x,j))\\
  &= \big(Q_M(F(R(x,j))) - \wt{T}(Q_M(F(x,j)))\big)\\ 
  &\ \ \ \ \ \ + \big(Q_N(F(R(x,j))) - \wt{T}(Q_N(F(x,j)))\big)\\
  &= \Big(\sum_{k=0}^\infty \wt{T}^k(Q_M(G(R^{-k}(x,j)))) 
      - \sum_{k=0}^\infty \wt{T}^{k+1}(Q_M(G(R^{-k-1}(x,j))))\Big)\\
  &\ \ \ \ \ \ + \Big(-\sum_{k=1}^\infty \wt{T}^{-k}(Q_N(G(R^k(x,j)))) 
      + \sum_{k=1}^\infty \wt{T}^{-k+1}(Q_N(G(R^{k-1}(x,j))))\Big)\\
  &= Q_M(G(x,j)) + Q_N(G(x,j))\\
  &= G(x,j),
\end{align*}
showing that $\Psi(F) = G$. Thus, $\Psi$ is bijective and its inverse is given by (\ref{eq5}).
Moreover, the estimates in (\ref{eq6}) and (\ref{eq7}) imply that
\begin{equation}\label{eq8}
\|\Psi^{-1}\| \leq \frac{a\,b\,(1+t)}{1-t}\,\cdot
\end{equation}

We will need maps of the form (\ref{eq4}) for two choices of $R$, namely: $R = \wt{T}$ and $R = S$.
In other words, we will consider the bijective continuous linear maps
\[
\Psi_1 : F \in \cF \mapsto F \circ \wt{T} - \wt{T} \circ F \in \cG \ \ \text{ and } \ \
\Psi_2 : F \in \cF \mapsto F \circ S - \wt{T} \circ F \in \cG.
\]
Let $\ov{\cL} \in \cG$ be defined by
\[
\ov{\cL}(x,j) = (L_j(x),j +_p 1).
\]
We define the maps $\Phi_1 : \cF \to \cG$ and $\Phi_2 : \cF \to \cG$ by
\[
\Phi_1(F)(x,j) = \ov{\cL}((x,j) + F(x,j)) \ \ \text{ and } \ \ \Phi_2(F)(x,j) = \Phi_1(F)(x,j) - \ov{\cL}(x,j).
\]
We claim that $\Phi_1$ and $\Phi_2$ are Lipschitz maps with
\begin{equation}\label{eq9}
\Lip(\Phi_1) = \Lip(\Phi_2) \leq \eps.
\end{equation}
Indeed, take $F,G \in \cF$ with $F(x,j) = (f_j(x),j)$ and $G(x,j) = (g_j(x),j)$. Since
\begin{align*}
|\ov{\cL}((x,j) + F(x,j)) - \ov{\cL}((&x,j) + G(x,j))|\\
&= |(L_j(x + f_j(x)), j +_p 1) - (L_j(x + g_j(x)),j +_p 1)|\\
&= \|L_j(x + f_j(x)) - L_j(x + g_j(x))\|\\
&\leq \Lip(L_j) \|f_j(x) - g_j(x)\|\\
&\leq \eps\, |F - G|_\infty,
\end{align*}
we obtain $\Lip(\Phi_1) \leq \eps$. It is clear that $\Lip(\Phi_2) = \Lip(\Phi_1)$.

By (\ref{eq8}) and (\ref{eq9}), the map $\Psi_1^{-1} \circ \Phi_1 : \cF \to \cF$ satisfies
\begin{equation}\label{eq10}
\Lip(\Psi_1^{-1} \circ \Phi_1) \leq \|\Psi_1^{-1}\| \Lip(\Phi_1) \leq \frac{a\,b\,(1+t)}{1-t}\, \eps \leq \delta < 1.
\end{equation}
Hence, by Banach's fixed point theorem, $\Psi_1^{-1} \circ \Phi_1$ has a unique fixed point $U$ in $\cF$.
Write $U(x,j) = (u_j(x),j)$ and let $H : \wt{X} \to \wt{X}$ be the continuous map defined by
\[
H(x,j) = (x + u_j(x),j).
\]
We claim that
\begin{equation}\label{eq11}
S \circ H = H \circ \wt{T}.
\end{equation}
Indeed, since $U$ is a fixed point of $\Psi_1^{-1} \circ \Phi_1$,
\[
U \circ \wt{T} - \wt{T} \circ U = \Psi_1(U) = \Phi_1(U).
\]
Hence,
\[
(u_{j +_p 1}(T(x)),j +_p 1) - (T(u_j(x)),j +_p 1) = (L_j(x + u_j(x)),j +_p 1),
\]
and so
\[
u_{j +_p 1}(T(x))  = T(u_j(x)) + L_j(x + u_j(x)).
\]
Therefore,
\begin{align*}
(S \circ H)(x,j) 
&= (S_j(x + u_j(x)),j +_p 1)\\
&= (T(x) + T(u_j(x)) + L_j(x + u_j(x)),j +_p 1)\\
&= (T(x) + u_{j +_p 1}(T(x)),j +_p 1)\\
&= (H \circ \wt{T})(x,j).
\end{align*}

Now, let $V = - \Psi_2^{-1}(\ov{\cL}) \in \cF$. 
Write $V(x,j) = (v_j(x),j)$ and let $K : \wt{X} \to \wt{X}$ be the continuous map defined by
\[
K(x,j) = (x + v_j(x),j).
\]
We claim that
\begin{equation}\label{eq12}
K \circ S = \wt{T} \circ K.
\end{equation}
Indeed, by the definition of $V$,
\[
V \circ S - \wt{T} \circ V = \Psi_2(V) = - \ov{\cL}.
\]
Hence,
\[
(v_{j +_p 1}(S_j(x)),j +_p 1) - (T(v_j(x)),j +_p 1) = (-L_j(x),j +_p 1),
\]
and so
\[
L_j(x) + v_{j +_p 1}(S_j(x))  = T(v_j(x)).
\]
Therefore,
\begin{align*}
(K \circ S)(x,j) 
&= (S_j(x) + v_{j +_p 1}(S_j(x)),j +_p 1)\\
&= (T(x) + L_j(x) + v_{j +_p 1}(S_j(x)),j +_p 1)\\
&= (T(x) + T(v_j(x)),j +_p 1)\\
&= (\wt{T} \circ K)(x,j).
\end{align*}

Let us now show that
\begin{equation}\label{eq13}
K(H(x,j)) = (x,j) \ \ \text{ for all } (x,j) \in \wt{X}.
\end{equation}
Let $W \in \cF$ be given by
\[
W(x,j) = (w_j(x),j), \ \text{ where } w_j(x) = u_j(x) + v_j(x + u_j(x)).
\]
Note that
\[
(K \circ H)(x,j) = (x + w_j(x),j).
\]
By (\ref{eq11}) and (\ref{eq12}), $K \circ H \circ \wt{T} = \wt{T} \circ K \circ H$, that is,
\[
(T(x) + w_{j +_p 1}(T(x)),j +_p 1) = (T(x) + T(w_j(x)),j +_p 1).
\]
Hence,
\begin{align*}
\Psi_1(W)(x,j) 
&= W(\wt{T}(x,j)) - \wt{T}(W(x,j))\\
&= (w_{j +_p 1}(T(x)),j +_p 1) - (T(w_j(x)),j +_p 1) = O_\cG(x,j).
\end{align*}
Thus, $W = O_\cF$ and (\ref{eq13}) holds.

Let us now prove that
\begin{equation}\label{eq14}
H(K(x,j)) = (x,j) \ \ \text{ for all } (x,j) \in \wt{X}.
\end{equation}
Let $Z \in \cF$ be given by
\[
Z(x,j) = (z_j(x),j), \ \text{ where } z_j(x) = v_j(x) + u_j(x + v_j(x)).
\]
Note that
\[
(H \circ K)(x,j) = (x + z_j(x),j).
\]
By (\ref{eq11}) and (\ref{eq12}), $H \circ K \circ S = S \circ H \circ K$, that is,
\[
(S_j(x) + z_{j +_p 1}(S_j(x)),j +_p 1) = (S_j(x + z_j(x)),j +_p 1).
\]
This is equivalent to
\[
(T(x) + L_j(x) + z_{j +_p 1}(S_j(x)),j +_p 1) = (T(x) + T(z_j(x)) + L_j(x + z_j(x)),j +_p 1),
\]
which gives
\[
z_{j +_p 1}(S_j(x)) - T(z_j(x)) = L_j(x + z_j(x)) - L_j(x).
\]
Hence,
\begin{align*}
\Psi_2(Z)(x,j) 
&= Z(S(x,j)) - \wt{T}(Z(x,j))\\
&= (z_{j +_p 1}(S_j(x)),j +_p 1) - (T(z_j(x)),j +_p 1)\\
&= (L_j(x + z_j(x)),j +_p 1) - (L_j(x),j +_p 1)\\
&= \ov{\cL}((x,j) + Z(x,j)) - \ov{\cL}(x,j)\\
&= \Phi_2(Z)(x,j).
\end{align*}
Thus, $Z$ is a fixed point of $\Psi_2^{-1} \circ \Phi_2$. As in (\ref{eq10}), the map $\Psi_2^{-1} \circ \Phi_2 : \cF \to \cF$ satisfies
\[
\Lip(\Psi_2^{-1} \circ \Phi_2) \leq \delta < 1.
\]
Hence, $\Psi_2^{-1} \circ \Phi_2$ has a unique fixed point in $\cF$.
Since $(\Psi_2^{-1} \circ \Phi_2)(O_\cF) = O_\cF$, we conclude that $Z = O_\cF$ and (\ref{eq14}) holds.

By (\ref{eq13}) and (\ref{eq14}), $H$ is a homeomorphism with inverse $K$. By definition,
\[
H(x,0) = (h(x),0), \ \ \text{ where } h(x) = x + u_0(x).
\]
It follows that $h : X \to X$ is a homeomorphism and
\[
H^{-1}(x,0) = (h^{-1}(x),0).
\]
By (\ref{eq11}), $S = H \circ \wt{T} \circ H^{-1}$, which gives
\[
S^p = H \circ \wt{T}^p \circ H^{-1}.
\]
In particular,
\[
((S_{p-1} \circ \cdots \circ S_0)(x),0) = S^p(x,0) = (H \circ \wt{T}^p \circ H^{-1})(x,0) = (h(T^p(h^{-1}(x))),0).
\]
Thus, $h$ establishes a conjugation between $S_{p-1} \circ \cdots \circ S_0$ and $T^p$. Finally,
\begin{equation}\label{extra}
\|h - I_X\|_\infty = \|u_0\|_\infty \leq |U|_\infty = |\Psi_1^{-1}(\Phi_1(U))|_\infty \leq \frac{a\,b\,(1+t)}{1-t}\, |\ov{\cL}|_\infty < \infty.
\end{equation}
This completes the proof of (A). 
Note that we have not assumed that the $L_j$'s are small in the sup norm so far.
However, this assumption will be necessary to establish (B), as the estimate obtained in (\ref{extra}) makes clear.
In fact, suppose that we also have
\[
\|L_j\|_\infty < \eps \text{ for all } j \in \{0,\ldots,p-1\}.
\]
Then, by (\ref{extra}),
\[
\|h - I_X\|_\infty < \frac{a\,b\,(1+t)}{1-t}\, \eps \leq \delta,
\]
which establishes (B).

Let us assume that $T$ is hyperbolic and prove the uniqueness of the homeomorphism~$h$.
Let $g : X \to X$ be any homeomorphism satisfying
\begin{equation}\label{eq15}
S_{p-1} \circ \cdots \circ S_0 = g \circ T^p \circ g^{-1} \ \text{ and } \ \ \|g - I_X\|_\infty < \infty.
\end{equation}
We define, recursively, maps $b_0,\ldots,b_{p-1}$ from $X$ into itself as follows:
\[
b_0 = g - I_X \ \ \text{ and } \ \ b_j = S_{j-1} \circ (I_X + b_{j-1}) \circ T^{-1} - I_X \ \text{ for } j \in \{1,\ldots,p-1\}.
\]
Clearly, these maps are continuous. Since $b_0$ is bounded and
\begin{equation}\label{eq16}
b_j(x) = T(b_{j-1}(T^{-1}(x))) + L_{j-1}(T^{-1}(x) + b_{j-1}(T^{-1}(x))),
\end{equation}
it follows recursively that all the $b_j$'s are bounded. Thus,
\[
B : (x,j) \in \wt{X} \mapsto (b_j(x),j) \in \wt{X}
\]
is a continuous bounded map. Since we are assuming that $T$ is hyperbolic, $Y = X$, and so $B \in \cF$. If we prove that
\begin{equation}\label{eq17}
B \text{ is a fixed point of } \Psi_1^{-1} \circ \Phi_1,
\end{equation}
then the uniqueness of such a fixed point would imply $B = U$ and we would obtain
\[
g = I_X + b_0 = I_X + u_0 = h,
\]
as desired. Given $(x,j) \in \wt{X}$, we have to show that
\[
\Psi_1(B)(x,j) = \Phi_1(B)(x,j).
\]
If $j \in \{0,\ldots,p-2\}$, then (\ref{eq16}) gives
\[
b_{j+1}(T(x)) = T(b_j(x)) + L_j(x + b_j(x)),
\]
and so
\begin{align*}
\Psi_1(B)(x,j)
&= B(\wt{T}(x,j)) - \wt{T}(B(x,j))\\
&= (b_{j+1}(T(x)) - T(b_j(x)),j + 1)\\
&= (L_j(x + b_j(x)),j + 1)\\
&= \Phi_1(B)(x,j).
\end{align*}
In order to analyze the case $j = p - 1$, first note that the recursive definition of the $b_j$'s implies that
\[
b_j = S_{j-1} \circ \cdots \circ S_0 \circ (I_X + b_0) \circ T^{-j} - I_X \ \text{ for all } j \in \{1,\ldots,p-1\}.
\]
In particular,
\[
I_X + b_{p-1} = S_{p-2} \circ \cdots \circ S_0 \circ g \circ T^{-p + 1}.
\]
Consequently,
\begin{align*}
T(b_{p-1}(x)) + L_{p-1}(x + b_{p-1}(x))
&= S_{p-1}(x + b_{p-1}(x)) - T(x)\\
&= (S_{p-1} \circ \cdots \circ S_0 \circ g)(T^{-p+1}(x)) - T(x)\\
&= (g \circ T^p)(T^{-p+1}(x)) - T(x)\\
&= g(T(x)) - T(x)\\
&= b_0(T(x)),
\end{align*}
where in the third equality we used the equality in (\ref{eq15}). Therefore,
\begin{align*}
\Psi_1(B)(x,p-1)
&= B(\wt{T}(x,p-1)) - \wt{T}(B(x,p-1))\\
&= (b_0(T(x)) - T(b_{p-1}(x)),0)\\
&= (L_{p-1}(x + b_{p-1}(x)),0)\\
&= \Phi_1(B)(x,p-1).
\end{align*}
This completes the proof of (\ref{eq17}), as it was to be shown.

Now, suppose that the operator $T$ is not hyperbolic. 
Then, we can take a nonzero vector $y \in M \cap T(N)$, and so
\[
z = \sum_{n = -\infty}^\infty T^n(y)
\]
defines a nontrivial fixed point of $T$. For each $\lambda \in \K$,
\[
h_\lambda : x \in X \mapsto h(x + \lambda z) \in X
\]
is a homeomorphism whose inverse is given by
\[
h_{\lambda}^{-1} : x \in X \mapsto h^{-1}(x) - \lambda z \in X.
\]
Since $z$ is a fixed point of $T$,
\[
h_\lambda \circ T^p \circ h_{\lambda}^{-1} = h \circ T^p \circ h^{-1} = S_{p-1} \circ \cdots \circ S_0 \ \ \text{ for all } \lambda \in \K.
\]
Moreover, 
\[
\|h_\lambda - I_X\|_\infty \leq \|h - I_X\|_\infty + |\lambda| \|z\| < \infty.
\]
If $\|h - I_X\|_\infty < \delta$, then we will also have $\|h_\lambda - I_X\|_\infty < \delta$ 
whenever $|\lambda|$ is sufficiently small.
In any case, we see that $h$ can be replaced by uncountably many other homeomorphisms.

Finally, let us prove the last assertion of the theorem.
For this purpose, we denote the maps $\Phi_1$, $U$ and $h$ by $\Phi_{\cL,1}$, $U_\cL$ and $h_\cL$, respectively,
in order to make it clear that they actually depend on $\cL$. Recall that:
\begin{itemize}
\item [(a)] $\Phi_{\cL,1} : \cF \to \cG$ is given by $\Phi_{\cL,1}(F)(x,j) = \ov{\cL}((x,j) + F(x,j))$.
\item [(b)] $\Psi_1^{-1} \circ \Phi_{\cL,1} : \cF \to \cF$ satisfies $\Lip(\Psi_1^{-1} \circ \Phi_{\cL,1}) \leq \delta < 1$.
\item [(c)] $U_\cL \in \cF$ is the unique fixed point of $\Psi_1^{-1} \circ \Phi_{\cL,1}$.
\item [(d)] $h_\cL : X \to X$ is the homeomorphism given by $h_\cL(x) = x + u_{\cL,0}(x)$,
  where $U_\cL(x,j) = (u_{\cL,j}(x),j)$.
\end{itemize}
Let $\Lambda$ denote the set of all $p$-tuples $\cL = (L_0,\ldots,L_{p-1})$ of bounded Lipschitz maps $L_j : X \to X$
with $\Lip(L_j) < \eps$ (and $\|L_j\|_\infty < \eps$ in the case of item (B)) endowed with the metric
\[
D(\cL,\cL') = \max_{0 \leq j < p} \|L_j - L'_j\|_\infty.
\]
We want to prove that the correspondence $\cL \to h_\cL$ is uniformly continuous in the sense described in the statement of the theorem.
Since
\begin{equation}\label{eq18}
\|h_\cL - h_{\cL'}\|_\infty = \|u_{\cL,0} - u_{\cL',0}\|_\infty \leq |U_\cL - U_{\cL'}|_\infty,
\end{equation}
it is enough to show that the map
\[
\cL \in \Lambda \mapsto U_\cL \in \cF
\]
is uniformly continuous.
In view of (b) and (c), this will follow from the parametrized version of Banach's fixed point theorem (see \cite{AGraJDug03}, for instance)
as soon as we show that the map
\begin{equation}\label{eq19}
f : (F,\cL) \in \cF \times \Lambda \mapsto \Psi_1^{-1}(\Phi_{\cL,1}(F)) \in \cF
\end{equation}
is uniformly continuous. But this is true because simple computations show that
\begin{equation}\label{eq20}
|f(F,\cL) - f(F',\cL')|_\infty \leq \|\Psi_1^{-1}\| (\eps |F - F'|_\infty + D(\cL,\cL')).
\end{equation}
This completes the proof of the theorem.
\end{proof}

\begin{remark}
Actually, the correspondence $\cL \to h_\cL$ is Lipschitz. In fact, by (\ref{eq19}) and~(c), we have that
\[
f(U_\cL,\cL) = \Psi_1^{-1}(\Phi_{\cL,1}(U_\cL)) = U_\cL \ \ \text{ for all } \cL \in \Lambda.
\]
In view of (b), it follows from the proof of the parametrized version of Banach's fixed point theorem that
\[
|U_{\cL'} - U_\cL|_\infty \leq \frac{2}{1 - \delta}\, |f(U_{\cL'},\cL') - f(U_{\cL'},\cL)|_\infty.
\]
Hence, by (\ref{eq8}), (\ref{eq18}) and (\ref{eq20}), we obtain
\begin{equation}\label{LipConstant}
\|h_\cL - h_{\cL'}\|_\infty \leq \frac{2\,a\,b\,(1+t)}{(1 - \delta)(1 - t)}\, D(\cL,\cL').
\end{equation}
It is worth noting that the Lipschitz constant in (\ref{LipConstant}) is independent of $p$.
\end{remark}

\begin{remark}
Theorem~\ref{Main2} remains true if we replace ``homeomorphism'' by ``uniform homeomorphism'' throughout its statement.
The proof is essentially the same, but we need to work with the function spaces $\cF$ and $\cG$ defined as follows:
\begin{itemize}
\item $\cF$ is the set of all uniformly continuous bounded maps $F : \wt{X} \to \wt{X}$ such that
  \[
  F(X \times \{j\}) \subset Y \times \{j\} \ \text{ for all } j \in \Z_p.
  \]
\item $\cG$ is the set of all uniformly continuous bounded maps $F : \wt{X} \to \wt{X}$ such that
  \[
  F(X \times \{j\}) \subset X \times \{j +_p 1\} \ \text{ for all } j \in \Z_p.
  \]
\end{itemize}
\end{remark}

\begin{remark}
Recall that a diffeomorphism $f : M \to M$ on a compact manifold is said to be {\em absolutely structurally stable} \cite{JFra73,JGuc72}
if there exist a neighborhood $\cN$ of $f$ in the $C^1$-topology and a constant $C \in (0,\infty)$ such that each $g \in \cN$ is
topologically conjugate to $f$ by means of a homeomorphism $h : M \to M$ satisfying
\[
\|h - I_X\|_\infty \leq C\, \|g - f\|_\infty.
\]
By using the notations in the proof of Theorem~\ref{Main2}, it follows from (\ref{extra}) that the constant
\[
C = \frac{a\,b\, (1+t)}{1-t}
\]
has the following property: for each $\cL = (L_0,\ldots,L_{p-1}) \in \Lambda$, the homeomorphism $h_\cL : X \to X$
conjugating $S_{p-1} \circ \cdots \circ S_0$ and $T^p$ (where $S_j = T + L_j$) satisfies
\[
\|h_\cL - I_X\|_\infty \leq C \max_{0 \leq j < p} \|L_j\|_\infty = C \max_{0 \leq j < p} \|S_j - T\|_\infty.
\]
\end{remark}

%%%%%%%%%%%%%%%%%%%%%%%%%%%%%%%%%%%%%%%%%%%%%%%%%%%%%%%%%%%%%%%
%%%%%%%%%%%%%%%%%%%%%%%%%%%%%%%%%%%%%%%%%%%%%%%%%%%%%%%%%%%%%%%

\section*{Declarations}

\noindent
{\bf Conflict of interest} \ The author declares that he has no conflict of interest.

%%%%%%%%%%%%%%%%%%%%%%%%%%%%%%%%%%%%%%%%%%%%%%%%%%%%%%%%%%%%%%%
%%%%%%%%%%%%%%%%%%%%%%%%%%%%%%%%%%%%%%%%%%%%%%%%%%%%%%%%%%%%%%%

\section*{Acknowledgement}

The author is beneficiary of a grant within the framework of the grants for the retraining, modality Mar\'ia Zambrano, 
in the Spanish university system (Spanish Ministry of Universities, financed by the European Union, NextGenerationEU).
The author was also partially supported by CNPq -- Project {\#}308238/2021-4, by CAPES -- Finance Code 001,
and by MCIN/AEI/10.13039/501100011033/FEDER, UE, Project PID2022-139449NB-I00.

%%%%%%%%%%%%%%%%%%%%%%%%%%%%%%%%%%%%%%%%%%%%%%%%%%%%%%%%%%%%%%%
%%%%%%%%%%%%%%%%%%%%%%%%%%%%%%%%%%%%%%%%%%%%%%%%%%%%%%%%%%%%%%%

\end{document}